\documentclass[12pt]{article}
\usepackage[latin1]{inputenc}
\usepackage{amsmath}
\usepackage{amssymb}
\usepackage{amsthm}

\numberwithin{equation}{section}

\begin{document}

\title{A New Proof of Stirling's Formula}
 \author{Thorsten Neuschel}
 
\date{}

\maketitle

\paragraph{Abstract} A new simple proof of Stirling's formula via the partial fraction expansion for the tangent function is presented.

\section{Introduction.}

Various proofs for Stirling's formula

\begin{equation}\label{1}n!\, \sim \, n^n\,e^{-n}\,\sqrt{2\,\pi\,n},~~~\text{as}~~n\rightarrow \infty,
\end{equation}
have been established in the literature since the days of de Moivre and Stirling in 1730 (for a historical exposition see, e.g., \cite{Diaconis}). Many of these proofs show that the limit 
\[\lim_{n\rightarrow\infty} \frac{n!}{n^n\,e^{-n}\,\sqrt{n}}\]
exists (for instance via the Euler-Maclaurin formula) in order to identify this limit by using the asymptotical behavior of the Wallis product, which is the crucial step. We will show that this last quite wily step can be replaced by a simple straightforward computation of the limit only using the partial fraction expansion for the tangent function 
\begin{equation}\label{2}
\pi \tan \pi x = \sum_{\nu=0}^\infty \frac{2x}{(\nu+\frac{1}{2})^2-x^2}\,.
\end{equation}
This expansion probably was found by Euler by the time Stirling determined his proof via Wallis' formula, see, e.g., \cite{Remmert}, p. 327.
For some alternative elementary proofs of Stirling's formula see, e.g., \cite{Diaconis}, \cite{Feller}, \cite{Michel}, \cite{Patin}, \cite{Robbins}.

\section{Proof.}

An application of the well-known Euler-Maclaurin formula in its simplest form (see, e.g., \cite{Wong}, p. 37, (6.21)) yields

\[\log n! = n\log n -n +1 + \log \sqrt{n} +\int\limits_{0}^{n-1} \frac{x-[x]-\frac{1}{2}}{1+x}\,dx.\]
In order to prove (\ref{1}), it is sufficient to show 

\begin{equation}\label{3}
\int\limits_{0}^{\infty} \frac{x-[x]-\frac{1}{2}}{1+x}\,dx = \log \sqrt{2\pi} -1.
\end{equation}
To prove this, we will show directly the identity

\begin{equation}\label{4}
\int\limits_{0}^{\infty} \frac{x-[x]-\frac{1}{2}}{1+x}\,dx = \int\limits_{0}^{1/2}\left( \frac{8x^2}{1-4x^2}-\pi x \tan \pi x\right)\,dx ,
\end{equation}
where the integral on the right-hand side can be evaluated by elementary calculus. We start our computation with

\begin{align*}
\int\limits_{0}^{\infty} \frac{x-[x]-\frac{1}{2}}{1+x}\,dx&=\sum_{\nu=0}^{\infty} \left\{\int\limits_{\nu}^{\nu+1/2} \frac{x-\nu-\frac{1}{2}}{1+x}\,dx + \int\limits_{\nu+1/2}^{\nu+1} \frac{x-\nu-\frac{1}{2}}{1+x}\,dx\right\}\\
&=\sum_{\nu=0}^{\infty} \left\{\int\limits_{0}^{1/2} \frac{x-\frac{1}{2}}{1+\nu+x}\,dx + \int\limits_{0}^{1/2} \frac{x}{\frac{3}{2}+\nu+x}\,dx\right\}.
\end{align*}
By an easy change of variables we observe
\[\int\limits_{0}^{1/2} \frac{x-\frac{1}{2}}{1+\nu+x}\,dx=-\int\limits_{0}^{1/2} \frac{x}{\frac{3}{2}+\nu-x}\,dx,\]
so that we obtain
\begin{align}\nonumber
\int\limits_{0}^{\infty} \frac{x-[x]-\frac{1}{2}}{1+x}\,dx&=\sum_{\nu=0}^{\infty}\int\limits_{0}^{1/2}\left(\frac{x}{\frac{3}{2}+\nu+x} -\frac{x}{\frac{3}{2}+\nu-x}\right)\,dx\\
\nonumber
&=\sum_{\nu=0}^{\infty}\int\limits_{0}^{1/2}\frac{-2x^2}{(\nu+\frac{3}{2})^2-x^2}\,dx\\
\label{5}
&=\int\limits_{0}^{1/2}\sum_{\nu=1}^{\infty}\frac{-2x^2}{(\nu+\frac{1}{2})^2-x^2}\,dx,
\end{align}
where the interchange of summation and integration is allowed due to the uniform convergence of the series in (\ref{5}) on the interval \([0,\frac{1}{2}]\). Applying (\ref{2}), we immediately obtain (\ref{4}). At this point of the proof, we have reduced the problem of determining the constant in Stirling's formula to a simple matter of elementary calculus as the resulting integral in (\ref{4}) can be evaluated easily. For convenience we will give some details. For example, using the decomposition
\[\frac{8x^2}{1-4x^2}=\frac{1}{1+2x}+\frac{1}{1-2x}-2\]
it can be rewritten as
\begin{align*}
\int\limits_{0}^{1/2}\left( \frac{8x^2}{1-4x^2}-\pi x \tan \pi x\right)\,dx &=\log \sqrt{2} -1 +\int\limits_{0}^{1/2}\left( \frac{1}{1-2x}-\pi x \tan \pi x\right)\,dx.
\end{align*}
Now, by a standard argumentation involving integration by parts, we can observe for \(0<\epsilon< 1/2\) that
\begin{align*}
\int\limits_{0}^{\epsilon}\left( \frac{1}{1-2x}-\pi x \tan \pi x\right)\,dx=&-\frac{1}{2} \log (1-2\epsilon)-\int\limits_{0}^{\epsilon}\pi x \tan \pi x\,dx\\
=& \left(\epsilon-\frac{1}{2}\right)\log \cos(\pi \epsilon)+\frac{1}{2}\log\frac{\cos(\pi \epsilon)}{1-2\epsilon}\\
&~~~-\int\limits_{0}^{\epsilon}\log \cos(\pi x)\,dx.
\end{align*}
Letting \(\epsilon\) tend to \(1/2\) we immediately obtain
\begin{align*}
\int\limits_{0}^{1/2}\left( \frac{1}{1-2x}-\pi x \tan \pi x\right)\,dx=\log \sqrt{\pi} -\log \sqrt{2} -\int\limits_{0}^{1/2}\log \cos (\pi x)\,dx.
\end{align*}
The remaining integral on the right-hand side can be evaluated elementary to \(-\log \sqrt{2}\) as shown, e.g., in \cite{Feller}, \cite{Feller2}. This computation relies on the fact that its value, say \(c\), remains unchanged if \(\cos (\pi x)\) is replaced by \(\sin (\pi x)\) so that we have (using the double angle formula)
\[\int\limits_{0}^{1/2}\log \sin(2 \pi x)\,dx=\log \sqrt{2}+2 \int\limits_{0}^{1/2}\log \sin(\pi x)\,dx.\]
As both integrals in the last equation coincide, we obtain \(c=-\log \sqrt{2}\), which completes the proof of (\ref{1}).

\bigskip

\noindent\textit{Department of Mathematics,
University of Trier, D-54286 Trier, Germany\\
neuschel@uni-trier.de}

\end{document}